\documentclass[11pt, a4paper, UKenglish]{article}

\usepackage[UKenglish]{babel}
\usepackage{amssymb, amsmath, textcomp, amsthm}
\usepackage{comment}
\usepackage{bbm,dsfont}
\usepackage{german}
\usepackage[pdftex]{graphicx}
\usepackage[utf8]{inputenc}
\usepackage{authblk}
\numberwithin{equation}{section}
\title{Hilbert transform along measurable vector fields constant on Lipschitz curves: $L^p$ boundedness}
\author{Shaoming Guo}
\date{}

\def\R{\mathbf{R}}
\def\N{\mathbf{N}}

\def\Z{\mathbf{Z}}

\def\lesim{\lesssim}

\def\begineq{\begin{equation}}
\def\endeq{\end{equation}}

\theoremstyle{plain}
\newtheorem{thm}{Theorem}[section]
\newtheorem{prop}[thm]{Proposition}
\newtheorem{lem}[thm]{Lemma}
\newtheorem{cor}[thm]{Corollary}
\newtheorem{defi}[thm]{Definition}
\newtheorem{claim}[thm]{Claim}
\newtheorem{rem}[thm]{Remark}

\newtheorem*{openproblem*}{Open Problem}

\begin{document}
\maketitle

\begin{abstract}
We prove the $L^p$ ($p>3/2$) boundedness of the directional Hilbert transform
in the plane relative to measurable vector fields which are constant
on suitable Lipschitz curves, extending the $L^2$ bounds in \cite{Guo}.

\end{abstract}

\let\thefootnote\relax\footnote{Date: \date{\today}}

\section{Statement of the main result}

In \cite{Guo} we proved that the Hilbert transforms along measurable vector fields which are constant on a suitable family of Lipschitz curves are bounded in $L^2$. The main goal of this paper is to generalize the above $L^2$ bounds to $L^p$ for $p$ other than 2 in the same setting.

\begin{thm}[Main Theorem]\label{theorem}
For vector fields $v:\R^2\to \R^2$ of the form $(1, u(h))$ where $h:\R^2\to \R$ is a Lipschitz function such that
\begineq\label{NN1.1}
\|\nabla h-(1, 0)\|_{\infty}\le \epsilon_0\ll 1,
\endeq
and $u:\R\to \R$ is a measurable function such that 
\begineq
\|u\|_{\infty}\le 1,
\endeq
the associated Hilbert transform, which is defined as
\begineq\label{NN1.3}
H_v f(x):=\int_{\R}f(x-tv(x))dt/t,
\endeq
is bounded in $L^p$ for all $p>3/2$.
\end{thm}

The above result is a Lipschitz perturbation of the following result by Bateman and Thiele in \cite{BT}, which is further based on Bateman \cite{Ba}, \cite{Ba2}, Lacey and Li \cite{LL1}, \cite{LL1}:
\begin{thm}(\cite{BT})\label{theorem1.2}
Let $v:\R^2\to \R^2$ be a one-variable vector field, i.e. vector field of the form 
\begineq
v(x_1, x_2)=(1, u(x_1)),
\endeq
for some measurable function $u$, then the associated Hilbert transform is bounded in $L^p$ for all $p>3/2$.
\end{thm}

In our Main Theorem, if we take $h(x_1, x_2)=x_1$, then the vector field becomes $(1, u(h(x_1, x_2)))=(1, u(x_1))$, which is a one-variable vector field. However, we have one more assumption that $\|u\|_{\infty}\le 1$. To recover the result in Theorem \ref{theorem1.2}, we just need to apply the following unisotropic scaling
\begineq
x_1\to x_1, x_2\to \lambda x_2,
\endeq
and a simple limiting argument.\\

As we state our main result as a Lipschitz perturbation of the one-variable vector fields, in the following, we will explain separately why the one-variable vector fields are interesting and why we do the perturbation at the level of the Lipschitz regularity but not others.

First of all, there is an interesting connection between the Hilbert transform along the one-variable vector fields and Carleson's maximal operator, which was observed by Coifman and El Kohen, we review the discussion as presented in \cite{BT}. Take a one-variable vector field $v(x_1, x_2)=(1, u(x_1))$, consider the associated Hilbert transform, which is given by
\begineq\label{NN1.6}
H_v f(x_1, x_2)=\int_{\R}f(x_1-t, x_2-t u(x_1))dt/t.
\endeq
Denoting by $\widehat{f}$ the partial Fourier transform in the second variable we obtain formally
\begin{equation}\label{onevarvf}
 \int f(x_1-t, x_2-u(x_1)t)\frac{dt}{t}
\end{equation}
$$ = \int e^{i x_2 \xi_2} \int \widehat{f}(x_1-t, \xi_2)
e^{iu(x_1)t\xi_2} \frac{dt}{t}d\xi_2.$$
By the Plancherel theorem, 
\begineq\label{NN1.8}
\|H_v f\|_2=\|\int \widehat{f}(x_1-t, \xi_2)
e^{iu(x_1)t\xi_2} \frac{dt}{t}\|_2
\endeq
For each fixed $\xi_2$, we recognize this to
essentially be the linearization of Carleson's maximal operator
\begin{equation}
(C f)(x):=\sup_{N\in \R}|\int_{\R}f(x-t)e^{iNt}\frac{dt}{t}|.
\end{equation} 
Hence the right hand side of \eqref{NN1.8} can be bounded by
\begineq\label{NN1.10}
\|C\hat{f}(x_1, \xi_2)\|_2\lesim \|\hat{f}(x_1, \xi_2)\|_2 \lesim \|f\|_2.
\endeq
Moreover, by choosing the function $u$ properly in \eqref{NN1.6}, the $L^2$ boundedness of $H_v$ also implies the $L^2$ boundedness of Carleson's maximal operator. 

Secondly, the class of the one-variable vector fields is also very natural from the viewpoint of the scaling symmetries. We leave the detailed discussion to the next section, where it will also become clear that the equivalence of  the $L^2$ bounds of $H_v$ and Carleson's maximal operator is due to the fact that they enjoy the same symmetries, especially the modulation symmetry.\\

Next we will explain the appearance of the Lipschitz regularity. For a vector field $v$, if one truncates \eqref{NN1.3} as 
\begineq\label{NN1.11}
H_{v, \epsilon_0}f(x):=\int_{-\epsilon_0}^{\epsilon_0}f(x-tv(x))dt/t,
\endeq
then it is reasonable to ask for pure regularity assumption on $v$ in order to bound $H_{v, \epsilon_0}$. Indeed, a counterexample in \cite{LL2}  based on the Perron tree construction of the Besicovitch-Kakeya set (see \cite{Pe} and \cite{Co}) shows that no bounds are possible for $v$ being H"older continuous of an exponent less than one, and it is a long standing open problem in harmonic analysis that whether Lipschitz regularity suffices.

At the regularity scale, the only known result is for analytic vector fields $v$ by Stein and Street in \cite{SS}, while the maximal variant of \eqref{NN1.11} in the same setting was proved much earlier by Bourgain \cite{Bo}. In the same direction, a prior result for smooth vector fields under certain geometric assumptions appeared in \cite{CNSW}. For some other partial results, see \cite{CSWW}, \cite{Katz2}, \cite{Ki}.

To our knowledge, the result in the present paper is the first that handles certain class of Lipschitz vector fields. Indeed, as has also been mentioned in \cite{Guo}, that our result has the following corollary, which includes a large class of Lipschitz vector fields:
\begin{cor}
For a measurable unit vector field $v_0:\R^2\to S^1$, suppose that 
\vspace{2mm}

i) there exists a bi-Lipschitz map $g_0: \R^2\to \R^2$ s.t.
\begineq
v_0(g_0(x_1, x_2)) \text{ is constant in } x_2;
\endeq 

ii) there exists $d_0>0$ s.t. $\forall x_1\in \R$,
\begineq
\angle (\partial_2 g_0(x_1, x_2), \pm v_0(g_0(x_1, x_2))) \ge d_0 \text{ for } x_2\text{-}a.e. \text{ in }\R.
\endeq 

\noindent Then the associated Hilbert transform is bounded in $L^p$ for all $p>3/2$, with the operator norm depending only on $p, d_0$ and the bi-Lipschitz norm of $g_0$. Moreover, the operator norm blows up when $d_0\to 0$.
\end{cor}

\begin{rem}\label{NNremark1.4}
The structure theorem for Lipschitz functions by Azzam and Schul in \cite{AS} states exactly that any Lipschitz function $u:\R^2\to \R$ (any Lipschitz unit vector field $v_0$ in our case) can be precomposed with a bi-Lipschitz function $g_0:\R^2\to \R^2$ such that $u\circ g_0$ is Lipschitz in the first coordinate and constant in the second coordinate, when restricted to a ``large'' portion of the domain.
\end{rem}

\vspace{4mm}

In the end, let us mention the new ingredients that will be used to extend the $L^2$ bounds in \cite{Guo}. Recall that in the $L^2$ case, the crucial ingredients are the use of Jones' beta numbers and the adapted $L^2$-Littlewood-Paley theory, which is in the spirit of the work on the Cauchy integral on Lipschitz curves (for example see \cite{CJS}). The techniques used in \cite{Guo} are the Hilbert space techniques as we need to use some facts like taking $L^2$ norm works trivially with certain square functions. Out of this reason, only $L^2$ bounds are obtained.

In the $L^p$ case for $p$ other than 2, one novelty is that we discovered a new paraproduct, which is indeed a one-parameter family of paraproducts, with each paraproduct living on one Lipschitz level curve of the vector field $v$. To prove the $L^p$ bounds for the one-parameter family of paraproducts, the difficulty is how to embed each paraproduct into two dimensions without losing orthogonality. To overcome this difficulty, we need to develop an adapted $L^p$-Littlewood-Paley theory, which again requires a new square function as an intermediate step. This new two dimensional square function shares some common features with the bi-parameter square function. See the following crucial Lemma \ref{lemma5.1} and Claim \ref{claim5.5}.

%
%
%
%

Another difference from the $L^2$ case in \cite{Guo} is that we will write the proof by using the $\delta$-calculus, which has been used intensively in the Fourier restriction estimates, see \cite{KF}, \cite{Fos} and \cite{CD} for example. One significant advantage of the $\delta$-calculus, which we will see shortly in the proof, is that it allows us to express everything in terms of the function $h$ from the Main Theorem, instead of going back and forth between $h$ and its inverse as in \cite{Guo}. For example, this can be seen by comparing the crucial definition of the adapted Littlewood-Paley operator associated to the vector fields, namely by comparing Definition 3.3 in \cite{Guo} with Definition \ref{defi3.3} in the current paper.\\


{\bf Organization of paper:} in Section 2, we will review the symmetries that were discussed in \cite{BT} for the Hilbert transforms along the one-variable vector fields. Moreover, we will introduce one more symmetry which appears only after we allow Lipschitz perturbation of the one-variable vector fields.

In Section 3 we will state the strategy of the proof for the Main Theorem.If we denote by $P_k$ a Littlewood-Paley operator in the second variable, the main observation in Bateman and Thiele's proof is that $H_v$ commutes with $P_k$. In our case, this is no longer true. To recover the orthogonality, an adapted Littlewood-Paley operator was introduced by the author in \cite{Guo} (see the following Definition \ref{defi3.3}), which allows to split the operator $H_v$ into a main term and a commutator term
\begineq
\sum_{k\in \Z} H_v P_k(f)=\sum_{k\in \Z}(H_v P_k(f)-\tilde{P}_k H_v P_k(f)+\tilde{P}_k H_v P_k(f)).
\endeq
The new symmetry is used in the definition of $\tilde{P}_k$.

The $L^p$ ($p>3/2$) bounds of the main term $\sum_{k\in \Z}\tilde{P}_k H_v P_k(f)$ can be proved essentially by the same argument as in Bateman and Thiele \cite{BT}, with just minor modifications that we will state in Section 4.

The main novelty is the $L^p$ boundedness of the commutator term 
\begineq\label{NN1.15}
\sum_{k\in \Z}(H_v P_k(f)-\tilde{P}_k H_v P_k(f)).
\endeq 
To achieve this, we will first review the time-frequency decomposition of the operator and the functions in Section 5, and then prove in Section 6 that \eqref{NN1.15} is bounded in $L^p$ for all $p>1$.\\

{\bf Acknowledgements.} The author would like to thank his advisor, Prof. Christoph Thiele, for his tremendous support. The author also thanks Diogo Oliveira e Silva for discussions on the $\delta$-calculus.

%
%
%
%
%

\section{Discussion on the symmetries}
In this section we will discuss various symmetries that the Hilbert transforms along vector fields have. We will start from the most general case, i.e. the case of measurable vector fields, and then introduce more and more assumptions suggested by the symmetries.

Given an arbitrary measurable vector field $v(x_1, x_2):\R^2\to \R^2$, by a renormalization, we assume that it is of the form $v(x_1, x_2)=(1, u(x_1, x_2))$ for some measurable function $u:\R^2\to \R$. Consider the associated Hilbert transform along this vector field, which is defined as
\begineq\label{KK2.3}
H_v f(x_1, x_2):=\int_{\R}f(x_1-t, x_2-t u(x_1, x_2))dt/t.
\endeq
Suppose for the moment that we would like to prove the following ideal estimate
\begineq
\|H_v f\|_p\le C \|f\|_p,
\endeq
for some $p>1$ and some universal constant $C$. We start by studying the symmetries of the above operators, which are, for example, translation, dilation and rotation. \\

First, it is simple to see that this operator is invariant under translation
\begineq
x_1\to x_1+x_{1,0}, x_2\to x_2+x_{2,0},
\endeq
with the vector field being changed to $(1, u(x_1-x_{1,0}, x_2-x_{2,0}))$, which is still a measurable function.

Next, we consider the dilation and rotation given by
\[ \left( \begin{array}{c}
x_1\\
x_2 \end{array}\right) 
= \left( \begin{array}{ccc}
a & b \\
d & e \end{array} \right)\cdot \left( \begin{array}{c}
x_1\\
x_2 \end{array}\right),\]
which are also supposed to be non-degenerate. By the decomposition of $2\times 2$ matrices, it is not difficult to see that there are in total four generators:
\[ A=\left( \begin{array}{cc}
\lambda_1 & 0 \\
0 & 1  \end{array}\right),  
B= \left( \begin{array}{ccc}
1 & 0 \\
0 & \lambda_2 \end{array} \right),
C=\left( \begin{array}{cc}
1 & 0\\
\lambda_3 & 1 \end{array}\right), 
D=\left( \begin{array}{cc}
1 & \lambda_4 \\
0 & 1 \end{array}\right).\]

{\bf Symmetry $\mathcal{A}$:} this is the dilation in the $x$ variable
\begineq
x_1\to \lambda_1 x_1, x_2\to x_2.
\endeq
Under this change of variables, the vector field is changed to $(1, \frac{1}{\lambda_1}u(\frac{x_1}{\lambda_1}, x_2))$.\\

{\bf Symmetry $\mathcal{B}$:} this is the dilation in the $y$ variable
\begineq\label{KK2.6}
x_1\to x_1, x_2\to \lambda_2 x_2.
\endeq
Under this change of variables, the vector field is changed to $(1, \frac{1}{\lambda_2}u(x_1, \frac{x_2}{\lambda_2}))$.\\

{\bf Symmetry $\mathcal{C}$:} this is what Bateman and Thiele called ``shearing transformation'' in \cite{BT}:
\begineq\label{symmetryc}
x_1\to x_1, x_2\to x_2+\lambda_3 x_1,
\endeq
with the vector field being changed to $(1, u(x_1, x_2-\lambda_3 x_1)+\lambda_3)$.\\

\begin{rem}
In frequency, the change of variables \eqref{symmetryc} corresponds to
\begineq
\xi_1\to \xi_1-\lambda_3 \xi_2, \xi_2\to \xi_2.
\endeq
Notice that if we restrict $\xi_2$ to one single frequency band, say $\xi_2\sim 1$, then roughly we have
\begineq
\xi_1\to \xi_1-\lambda_3,
\endeq
which is the translation in the frequency variable $\xi_1$. Indeed, it will become clear later in the time-frequency decomposition in Section \ref{subsection4.1} that this is the same as the modulation invariance in Carleson's maximal operator.
\end{rem}

{\bf Symmetry $\mathcal{D}$:} this is the shearing transformation with $x_1$ and $x_2$ being exchanged:
\begineq
x_1\to x_1+\lambda_4 x_2, x_2\to x_2,
\endeq
with the vector field being changed to $(1, u(x_1-\lambda_4 x_2, x_2)+\lambda_4)$.\\

So far we have shown that if we only assume the vector field to be measurable, then the operator \eqref{KK2.3} satisfies the translation symmetry and the Symmetries $\mathcal{A}$, $\mathcal{B}$, $\mathcal{C}$ and $\mathcal{D}$. Unfortunately, even for H"older continuous vector field (with exponent less than one), the operator \eqref{KK2.3} might not be bounded in $L^p$ for any $p\ge 1$.

However, if we eliminate the Symmetry $\mathcal{D}$ (or equivalently the Symmetry $\mathcal{C}$) from the class of the measurable vector fields, then it is not difficult to see that a very natural choice is the class of the one-variable measurable vector fields, which enjoys all the other symmetries. Moreover, as it has been pointed out before, that by exploring the translation symmetry and the symmetries $\mathcal{A}$, $\mathcal{B}$ and $\mathcal{C}$, Bateman \cite{Ba2}, Bateman and Thiele \cite{BT} have proved that the operator \eqref{KK2.3} is bounded in $L^p$ ($\forall p>3/2$) for arbitrary measurable one-variable vector fields.\\

Let us explain a bit more from the viewpoint of symmetries why we expect the operator in \eqref{KK2.3} to be bounded for the one-variable vector fields, but not for arbitrary measurable vector fields: for a one-variable vector $v(x_1, x_2)=(1, u(x_1))$, if we denote by $P_k$ the Littlewood-Paley projection operator in the vertical variable, then what Bateman has proved in \cite{Ba2} is 
\begineq
\|H_v P_k f\|_p\lesim \|P_k f\|_p, \forall p\in (1, \infty),
\endeq
with the bound being independent of $k\in \Z$. Notice that the operator $H_v P_k$ has the translation symmetry, symmetry $\mathcal{A}$ and $\mathcal{C}$, which correspond to the translation, dilation and modulatioin symmetries for Carleson's maximal operator. In this sense, we say that Bateman's result is equivalent with the boundedness of Carleson's maximal operator (the precise calculation is done in \eqref{onevarvf}-\eqref{NN1.10}). 

When trying to put all the frequency annuli together to prove the boundedness of the whole operator $H_v$, we have increased the complexity of the problem ``by one dimension''. Fortunately, this can be compensated by making use of the Symmetry $\mathcal{B}$, which is done by Bateman and Thiele in \cite{BT} through a square function estimate.

So far we have seen that for the one-variable vector fields, we have made use of the exact number of symmetries. However, for the arbitrary measurable vector fields, there is an extra Symmetry $\mathcal{D}$ that we should respect, which serves as the heuristic for not expecting \eqref{KK2.3} to be bounded.\\

However, we still want to go beyond the one-variable vector fields. Notice that for the Hilbert transform along a general Lipschitz vector field, both Symmetry $\mathcal{C}$ and Symmetry $\mathcal{D}$ might still appear at the same time: for $v(x)=(1, u(x_1, x_2))$ with $u$ being Lipschitz, by applying Symmetry $\mathcal{C}$ with $\lambda_3$ being small and Symmetry $\mathcal{D}$ with $\lambda_4$ being small, what we get is 
\begineq
(1, u(x_1-\lambda_3 x_2, x_2-\lambda_4 x_1-\lambda_3 \lambda_4 x_2)+\lambda_3+\lambda_4),
\endeq
which is still a Lipschitz vector field with a comparable Lipschitz constant. 

Indeed, by including a Lipschitz perturbation of the one-variable vector fields (see the assumption of the Main Theorem), we bring the Symmetry $\mathcal{D}$ into the problem, with the cost that all the symmetries $\mathcal{A}$-$\mathcal{D}$ become ``quasi-symmetries''. Let us explain what we mean by this: the vector fields that we can handle are of the form $(1, u(h(x)))$, where $\|u\|_{\infty}\le 1$ and $h$ is a Lipschitz function satisfying
\begineq
\|\nabla h-(1, 0)\|_{\infty}\le \epsilon_0\ll 1.
\endeq
If we apply the Symmetry $\mathcal{D}$ with 
\begineq\label{KK2.12}
\lambda_4\ll 1,
\endeq 
the new vector field $(1, u_{\lambda_4}(h_{\lambda_4}(x)))$ will satisfy
\begineq
\|u_{\lambda_4}\|\le 2, \|\nabla h_{\lambda_4}-(1, 0)\|_{\infty}\le 2\epsilon_0,
\endeq
i.e. under the action of the symmetry, the assumption on the vector field is preserved up to a factor of two. This explains the notion of ``quasi-symmetry''.

The new quasi-symmetry $\mathcal{D}$ will be used implicitly in Definition \ref{defi2.2}, hence it will also be used in the crucial Definition \ref{defi3.3} of the adapted Littlewood-Paley projection operators.

\section{Strategy of the proof of the Main Theorem}\label{section3}

We first observe that if we denote by $\Gamma$ the two-ended cone which forms an angle less than $\pi/4$ with the vertical axis, then by the assumption that $|u|\le 1$, we can w.l.o.g. assume that 
\begineq
\text{supp } \hat{f}\subset \Gamma,
\endeq
as for functions $f$ with frequency supported on $\R^2\setminus \Gamma$, we have that
\begineq
H_vf(x)=H_{(1, 0)}f(x),
\endeq 
which is the Hilbert transform along the constant vector field $(1, 0)$. But $H_{(1,0)}$ is bounded by Fubini's theorem and the $L^2$ boundedness of the Hilbert transform.\\

The rest of the proof consists of two relatively independent steps. The first step will just be an adaption of Bateman and Thiele's argument in \cite{BT} to our case. Our key observation is that both covering lemmas used there (Lemma 7 and Lemma 8) indeed hold true in our setting, from which we can derive the following proposition as a corollary by repeating the rest of the argument in \cite{BT}.
\begin{prop}\label{prop2.1}
Under the same assumptions as in the Main Theorem, we have the following square function estimate
\begineq\label{NN2.3}
\left\|\left(\sum_{k\in \Z}(H_v P_k(f))^2\right)^{1/2}\right\|_p\lesim \|f\|_p, \forall p>3/2,
\endeq
where $P_k$ is the $k$-th Littlewood-Paley projection operator in the vertical direction.
\end{prop}

\begin{rem}
The operator $P_k$ is defined in the following way: if we denote by $\psi_0$ is a smooth function with support on $[-5/2, -1/2]\cup [1/2, 5/2]$ such that 
\begineq
\sum_{k\in \Z}\psi_k(t)=1, \forall t\neq 0,
\endeq
and 
\begineq\label{NN2.5}
\psi_k(t):=\psi_0(2^{-k} t),
\endeq
then 
\begineq
P_k f(x_1, x_2):=\int_{\R}f(x_1, x_2-y_2)\check{\psi}_k(y_2)dy_2.
\endeq
\end{rem}

For the one-variable vector fields, i.e. vector fields of the form $v(x, y)=(1, u(x))$ for some measurable function $u$, Bateman and Thiele in \cite{BT} used \eqref{NN2.3} and the crucial observation that
\begineq\label{NN2.4}
H_v P_k=P_k H_v
\endeq
to conclude the boundedness of $H_v$. In our case, the identity \eqref{NN2.4} is no longer true, i.e. the orthogonality between $H_v P_k f$ for different $k\in \Z$ is missing. 

To recover the orthogonality, an adapted Littlewood-Paley operator along the level curves of the vector field was introduced by the author in \cite{Guo}. This operator is in the spirit of prior work on the Cauchy integral on Lipschitz curves, but more of a bi-parameter type as we have one-parameter family of level curves.

Here we give an equivalent definition of the operator $\tilde{P}_k$ by using the language of the $\delta$-calculus. The advantage of this new definition is, compared with the one in \cite{Guo}, that it does not necessitate neither the change of coordinates nor the parametrization of the Lipschitz curves, both of which can be replaced by introducing the following auxiliary function. To do this, we need several notations: for $t\in \R$ we define 
\begineq
\Gamma_t:=\{x\in \R^2: h(x)=t\}.
\endeq 
Moreover, we denote by $v_t$ the value of the vector field $v$, which is a constant along $\Gamma_t$.

\begin{defi}[Auxiliary Function]\label{defi2.2}
For every $t\in R$, we define a new function $h_t:\R^2\to \R$ in such a way that, if for some $y\in \Gamma_t$, we have 
\begineq
z-y=d\cdot v_t
\endeq 
for some $d\in \R$, then we set $h_t(z)=d.$
\end{defi}
\begin{rem}
It is not difficult to see that
\begineq\label{GG3.5}
|\nabla h_t|\sim 1, \text{ a.e. in }  \R^2,
\endeq
where the constant is independent of $t\in \R$.
\end{rem}

\begin{defi}[Adapted Littlewood-Paley Operator]\label{defi3.3}
For $x\in \R^2$, we denote $t=h(x)$. We then define the adapted Littlewood-Paley projection operator $\tilde{P}_k$ restricted on the curve $\Gamma_t$ by
\begineq\label{GG3.4}
\tilde{P}_k f(x) :=\int_{\R^2} \delta(h_t(y))f(y)\check{\psi}_k((x-y)\cdot v_t^{\perp})dy,
\endeq
where $\psi_k(\cdot)$ is given by \eqref{NN2.5}.
\end{defi}

\begin{rem}
We show that the above Definition \ref{defi3.3} is equivalent with the Definition 3.3 in \cite{Guo}. To do this, we start from the new definition \eqref{GG3.4}: for a fixed $t\in \R$, the two vectors $v_t$ and $v_t^{\perp}$ form a orthogonal coordinate system of the plane. Write $y\in \R^2$ in this new system as
\begineq
y=y_1 v_t +y_2 v_t^{\perp},
\endeq
and for the sake of simplicity we will still use the notation $y=(y_1, y_2)$. This changes the expression in \eqref{GG3.4} to 
\begineq\label{GG3.7}
\begin{split}
& \int_{\R^2} f(y_1, y_2)\delta(h_t(y_1, y_2))\check{\psi}_k(x_2-y_2)dy\\
& =\int_{\R}\left(\int_{\R} f(y_1, y_2)\delta(h_t(y_1, y_2))dy_1\right)\check{\psi}_k(x_2-y_2)dy_2.
\end{split}
\endeq
Hence if we use the same parametrization as the one in Definition 3.3 in \cite{Guo}, which is 
\begineq
\Gamma_t=\{y_2 v_t^{\perp}+ g_t(y_2)v_t |  y_2\in \R\},
\endeq
then by the definition of the function $h_t$ in Definition \ref{defi3.3}, which implies 
\begineq
\int_{\R}\delta(h_t(x))dx_1=1,
\endeq 
the right hand side of \eqref{GG3.7} will equal 
\begineq
\int_{\R}f(g_t(y_2), y_2)\check{\psi}_k(x_2-y_2)dy_2,
\endeq
which is exactly the one given by the Definition 3.3 in \cite{Guo}.
\end{rem}

\begin{lem}[Adapted Littlewood-Paley Theory]\label{littlewoodpaley}
For $p\in (1, \infty)$, we have the following variants of the Littlewood-Paley theorem:
\begin{align}
&\label{littlewood1}\|(\sum_{k\in \Z}|\tilde{P}_k f|^2)^{1/2}\|_p\sim \|f\|_p,\\
&\label{littlewood2} \|(\sum_{k\in \Z}|\tilde{P}^{*}_k f|^2)^{1/2}\|_p\sim \|f\|_p.
\end{align}
\end{lem}

\noindent {\bf Proof of Lemma \ref{littlewoodpaley}:} by the Fubini theorem, we obtain
\begineq
\int_{\R^2} \left(\sum_{k\in \Z}|\tilde{P}_k f|^2\right)^{p/2}=\int_{\R}\int_{\R^2}\left(\sum_{k\in \Z}|\tilde{P}_k f|^2\right)^{p/2} \delta(h(x)-t)dxdt.
\endeq
When integrating against $dx$, by doing the change of variables $h(x)-t\to h_t(x)$, we can write the right hand side of the above expression as
\begineq
\int_{\R}\int_{\R^2}\left(\sum_{k\in \Z}|\tilde{P}_k f|^2\right)^{p/2} \delta(h_t(x))\frac{|\nabla h_t(x)|}{|\nabla h(x)|}dxdt.
\endeq
By the bound on $\nabla h_t$ in \eqref{GG3.5} and our assumption on $\nabla h$ in \eqref{NN1.1} that 
\begineq
|\nabla h|\sim 1, \text{ a.e. in }  \R^2,
\endeq
it suffices to show that
\begineq\label{GG3.18}
\int_{\R^2} \left(\sum_{k\in \Z}|\tilde{P}_k f|^2\right)^{p/2}\delta(h_t(x)) dx \lesim \int_{\R^2}|f(x)|^p\delta(h_t(x))dx,
\endeq
with a bound being independent of $t\in \R$.

We substitute the definition of $\tilde{P}_k$ into the left hand side of the last expression to obtain
\begineq\label{EE2.9}
\int_{\R^2} \left(\sum_{k\in \Z} \left|\int_{\R^2} \delta(h_t(y))f(y)\check{\psi}_k((x-y)\cdot v_t^{\perp})dy\right|^2\right)^{p/2} \delta(h_t(x))dx.
\endeq
The above expression can be viewed as a two dimensional Littlewood-Paley operator with the singular measure $\delta(h_t(\cdot))$, hence heuristically it is bounded by 
\begineq
\int_{\R^2}|f(x)|^{p} \delta(h_t(x))dx.
\endeq

To make the above argument rigorous, we introduce the change of variables 
\begineq
x\to x_1 v_t+x_2 v_t^{\perp}, y\to y_1 v_t + y_2 v_t^{\perp}.
\endeq
For the sake of simplicity, after the change of variables, we will still write $x=(x_1, x_2)$ and $y=(y_1, y_2)$. The expression in \eqref{EE2.9} hence becomes
\begineq
\begin{split}
& \int_{\R^2} \left(\sum_{k\in \Z} \left|\int_{\R^2} \delta(h_t(y))f(y)\check{\psi}_k(x_2-y_2)dy\right|^2\right)^{p/2} \delta(h_t(x))dx\\
& =\int_{\R^2} \left(\sum_{k\in \Z} \left|\int_{\R}\left(\int_{\R} \delta(h_t(y))f(y)dy_1\right) \check{\psi}_k(x_2-y_2)dy_2\right|^2\right)^{p/2} \delta(h_t(x))dx.
\end{split}
\endeq
Notice that for any $x_2\in \R$, we have 
\begineq
\int_{\R}\delta(h_t(x))dx_1=1.
\endeq 
Hence the right hand side of the last display becomes
\begineq
\int_{\R} \left(\sum_{k\in \Z} \left|\int_{\R}\left(\int_{\R} \delta(h_t(y))f(y)dy_1\right) \check{\psi}_k(x_2-y_2)dy_2\right|^2\right)^{p/2}dx_2.
\endeq
It is not difficult to see that the above is just a one-dimensional Littlewood-Paley square function for the function 
\begineq
\int_{\R} \delta(h_t(y))f(y)dy_1,
\endeq
hence it can be bounded by
\begineq
\int_{\R}\left|\int_{\R} \delta(h_t(y))f(y)dy_1\right|^p dy_2=\int_{\R^2} \delta(h_t(x))|f(x)|^{p}dx.
\endeq
So far we have finished the proof of \eqref{GG3.18}, thus \eqref{littlewood1}. For the second equivalence relation \eqref{littlewood2}, the proof is similar, hence we leave it out. $\Box$\\

To proceed, we will split the operator into two terms, 
\begineq
\sum_{k\in \Z} H_v P_k(f)=\sum_{k\in \Z}(H_v P_k(f)-\tilde{P}_k H_v P_k(f)+\tilde{P}_k H_v P_k(f)).
\endeq
Then by the triangle inequality, we have
\begineq
\|\sum_{k\in \Z} H_v P_k(f)\|_p\lesim\|\sum_{k\in \Z}(H_v P_k(f)-\tilde{P}_k H_v P_k(f))\|_p+\|\sum_{k\in\Z}\tilde{P}_k H_v P_k(f)\|_p.
\endeq
We call the second term the main term, and the first term the commutator term.

To bound the main term, we first use duality to write the $L^p$ norm into
\begin{align*}
\|\sum_{k\in \Z} \tilde{P}_k H_v P_k(f)\|_p &= \sup_{\|g\|_{p^{\prime}}=1}|\langle \sum_{k\in \Z}\tilde{P}_k H_v P_k(f), g\rangle|\\
		&= \sup_{\|g\|_{p^{\prime}}=1}|\sum_{k\in \Z}\langle H_v P_k(f), \tilde{P}_k^{*}(g)\rangle|.
\end{align*}
Then by Cauchy-Schwartz and H"older's inequality, we bound the right hand side by
\begineq\label{NN2.31}
\begin{split}
& \sup_{\|g\|_{p^{\prime}}=1}\int (\sum_{k\in \Z}|H_v P_k(f)|^2)^{1/2}(\sum_{k\in \Z} |\tilde{P}_k^{*}(g)|^2)^{1/2}\\
& \lesim \sup_{\|g\|_{p^{\prime}}=1} \|(\sum_{k\in \Z}|H_v P_k(f)|^2)^{1/2}\|_p \|(\sum_{k\in \Z} |\tilde{P}_k^{*}(g)|^2)^{1/2}\|_{p^{\prime}}.
\end{split}
\endeq
In the end, by applying Proposition \ref{prop2.1} to the former term in the last expression and Lemma \ref{littlewoodpaley} to the latter term, we get the desired bound
\begineq
\eqref{NN2.31}\lesim \sup_{\|g\|_{p^{\prime}}=1} \|f\|_p \|g\|_{p'}=\|f\|_p.
\endeq

\vspace{4mm}

%
%
%
%
%

Now we turn to the commutator term. Before explaining the idea of estimating the commutator term, we recall some notations from \cite{Guo}. Select a Schwartz function $\psi_0$ such that $\psi_0$ is supported on $[\frac{1}{2},\frac{5}{2}]$, let 
\begineq\label{EE2.6}
\psi_l(t):=\psi_0(2^{-l} t).
\endeq
By choosing $\psi_0$ properly, we can construct a partition of unity for $\R^+$, i.e. 
\begineq
\mathbbm{1}_{(0,\infty)}=\sum_{l\in \Z} \psi_l.
\endeq
Let 
\begineq\label{E4.24}
H_{l}f(x):= \int \check{\psi}_l(t)f(x-tv(x))dt.
\endeq
Then the operator $H_v$ can be decomposed into the sum
\begineq\label{EE1.6}
H_v=-\mathds{1}+ 2 \sum_{l\in \Z} H_{l}.
\endeq

We continue to explain the strategy of proving the $L^p$ boundedness of the commutator term, which is
\begineq
\|\sum_{k\in \Z}(H_v P_k(f)-\tilde{P}_k H_v P_k(f))\|_p \lesim \|f\|_p.
\endeq
By the dyadic decomposition in \eqref{EE1.6}, this is equivalent to bound the following
\begineq
\sum_{k\in \Z} \sum_{l\in \Z}(H_{l} P_k f- \tilde{P}_k H_{l} P_k f).
\endeq
Notice that by definition, $H_{l} P_k f$ vanishes for $l>k$, which simplifies the last expression to 
\begineq
\sum_{l\ge 0} \sum_{k\in \Z} (H_{k-l} P_k f- \tilde{P}_k H_{k-l} P_k f).
\endeq
So by the triangle inequality it suffices to prove

\begin{prop}\label{prop2.4}
Under the same assumptions as in the Main Theorem, for any $p\in (1, \infty)$, there exists a constant $\gamma_p>0$ such that
\begineq\label{EE2.17}
\|\sum_{k\in \Z} (H_{ k-l} P_k(f)-\tilde{P}_k H_{k-l} P_k(f))\|_p \lesim 2^{-\gamma_p l}\|f\|_p,
\endeq
with the constant being independent of $l\in \N$.
\end{prop}

The idea of proving endpoint estimates like  the $L^{\infty}\to BMO$ estimate will probably not work as the output of the operator $H_v$ is so rough that it is only measurable across the family of Lipschitz level curves, in another word, the orthogonality between different tiles is missing.

To recover the orthogonality at the level of the $L^2$ estimate, the argument in \cite{Guo} relies heavily on the fact that taking $L^2$ norm works perfectly (also trivially) with the square function. Hence we could expand certain square summation and apply H"older's inequality to turn the problem to the analysis on every single Lipschitz curve. 

However, in the $L^p$ estimate for $p\neq 2$, this strategy does not work, and instead we will invoke a new square function as an intermediate step. This square function is similar to the square function in the product space $\R\times \R$. 

\begin{rem}
Although the endpoint $L^{\infty}\to BMO$ estimate might not work for \eqref{EE2.17} with the classical $BMO$ space,  we still hope that there would be some variants, possibly similar to the fiber-wise Hardy and $BMO$ spaces in \cite{Bernicot} and \cite{Kovac}, which will act as the right substitutes for the endpoint theory.
\end{rem}

\begin{rem}
For the one-variable vector fields $v(x_1, x_2)=(1, u(x_2))$, it was proved in \cite{CSWW}, under some convexity and curvature assumptions on the function $u:\R\to \R$, that the associated Hilbert transform and maximal function map $H^1_{prod}(\R\times \R)$ to $L^1$, where $H^1_{prod}(\R\times \R)$ denotes the product Hardy space.

However, it was also pointed out that this might not be the right endpoint theory, and some new underlying Calderon-Zygmund theory is to be expected. See Remark $(iii)$ in Page 597 in \cite{CSWW}.
\end{rem}

\section{Boundedness of the main term: proof of Proposition \ref{prop2.1}}\label{section4}

The goal of this section is to make an observation that Bateman and Thiele's square function estimate (see $(2.1)$ in \cite{BT}) for the one-variable vector fields, which is
\begineq\label{E3.1}
\|(\sum_{k\in \Z}(H_v P_k(f))^2)^{1/2}\|_p\lesim \|f\|_p, \forall p>3/2,
\endeq
 works equally well for our case, with just minor modifications. Indeed, the proof of the estimate \eqref{E3.1} is reduced by Bateman and Thiele in \cite{BT} to three covering lemmas (Lemma 7 and Lemma 8 in \cite{BT}, Lemma 6.2 in \cite{Ba2}), and our observation is that all these covering lemmas still hold true for the case where the vector fields are constant only on Lipschitz curves instead of vertical lines.

Before stating the covering lemmas and the modification that we will make in the proof, we first need to introduce several definitions.

\begin{defi}
For a rectangle $R\subset \R^2$, with $l_R$ its length, $w_R$ its width, we define its uncertainty interval $EX(R)\subset \R$ to be the interval of width $w_R/l_R$ and centered at slope($R$). Denote by $E(R)$ the collection of the points $x\in R$ s.t. the vector $v(x)=(1, u(h(x)))$ points roughly in the same direction as the long side of $R$:
\begineq
E(R)=\{x\in R: u(h(x))\in EX(R)\}.
\endeq
Then the popularity of the rectangle $R$ is defined to be 
\begineq
pop_R:=|\{x\in \R^2: u(h(x))\in EX(R)\}|/|R|.
\endeq
Here $u$ and $h$ are the two functions in the Main Theorem.
\end{defi}
\begin{defi}
Given two rectangles $R_1$ and $R_2$ in $\R^2$, we write $R_1\le R_2$ whenever $R_1\subset C R_2$ and $EX(R_2)\subset EX(R_1)$, where $C$ is some properly chosen large constant, and $C R_2$ is the rectangle with the same center as $R_2$ but dilated by the factor $C$. 
\end{defi}

Now we are ready to state the key covering lemmas:
\begin{lem}(Lemma 6.2 in \cite{Ba2}, see also Lemma 4.3 in \cite{Guo})
Suppose $\mathcal{R}_0$ is a collection of pairwise incomparable (under ``$ \le $'') rectangles of uniform width such that for each $R\in \mathcal{R}_0$,  we have
\begineq
pop_{R}\ge \delta \text{ and }\frac{1}{|R|}\int_R \mathbbm{1}_F \ge \lambda.
\endeq
Then under the same assumptions on $u$ and $h$ as in the Main Theorem, we have for each $p>1$ that 
\begineq
\sum_{R\in \mathcal{R}_0}|R| \lesim \frac{|F|}{\delta \lambda^p}.
\endeq
\end{lem}

\begin{lem}(Lemma 7 in \cite{BT})\label{lem3.2}
Under the same assumptions as in the Main Theorem, let $\delta>0$ and $q>1$, let $G\subset \R^2$ be a measurable set and $\mathcal{R}$ be a finite collection of rectangles such that 
\begineq\label{E3.3}
|E(R)\cap G|\ge \delta |G|
\endeq
for each $R\in \mathcal{R}$. Then
\begineq
|\bigcup_{R\in \mathcal{R}} R|\lesim \delta^{-q}|G|.
\endeq
\end{lem}

\begin{lem}(Lemma 8 in \cite{BT})
Under the same assumptions as in the Main Theorem, let $0<\sigma, \delta\le 1$, let $H$ be a measurable set, and let $\mathcal{R}$ be a finite collection of rectangles such that for each $R\in \mathcal{R}$ we have
\begineq
pop_R\ge \sigma, |H\cap R|\ge \delta |R|.
\endeq
Then
\begineq
|\bigcup_{R\in \mathcal{R}}R|\lesim \sigma^{-1} \delta^{-2} |H|.
\endeq
\end{lem}

To prove these covering lemmas, one just need to replace the classical rectangles by the following ``rectangles'' adapted to the vector fields, and run the same argument as in Bateman and Thiele in \cite{BT}. 
\begin{defi}(rectangles adapted to the vector field)\label{definition3.6}
For a rectangle $R\subset \R^2$, with its two long sides lying on the parallel lines $x_2=kx_1+b_1$ and $x_2=kx_1+b_2$ for some $k\in [-1, 1]$ and $b_1, b_2\in \R$, define $\tilde{R}$ to be the adapted version of $R$, which is given by the set
\begineq
\{x\in \R^2: h(x)\in h(R)\}\bigcap \{(x_1, kx_1+b):x_1\in \R, b\in [b_1, b_2]\},
\endeq
where $h:\R^2\to \R$ is the function from the Main Theorem.
\end{defi}

From $R$ to $\tilde{R}$, the  length and the width of the rectangle are preserved up to a constant, and the same also holds true for the ``popularity''. Moreover, the proofs in \cite{BT} are ``stable'' under bi-Lipschitz mapping. Hence we will leave out the details and refer to \cite{BT}.\\

These two lemmas were used to give an upper bound on the size of the exceptional sets around which the rectangles have either large size or large density. After excluding the exceptional sets, the argument in \cite{BT}, together with \cite{Ba2}(which also works equally well for our case as has been pointed out in \cite{Guo}), will lead to the square function estimate, i.e. Proposition \ref{prop2.1}.

\section{Time-frequency decomposition}\label{subsection4.1}
The content of this section is the same as the Subsection 5.1 in \cite{Guo}. We still include these notations here for the sake of completeness.\\

{\bf Discretizing the functions:} Fix $l\ge 0$, we write $\mathcal{D}_l$ as the collection of the dyadic intervals of length $2^{-l}$ contained in $[-2,2]$. Fix a smooth positive function $\beta:\R\to \R$ s.t. 
\begineq
\beta(x)=1, \forall |x|\le 1; \beta(x)=0, \forall |x|\ge 2.
\endeq
Also choose $\beta$ such that $\sqrt{\beta}$ is a smooth function. Then fix an integer $c$(whose exact value is unimportant), for each $\omega\in \mathcal{D}_l$, define
\begineq
\beta_{\omega}(x)= \beta(2^{l+c}(x-c_{\omega_1})),
\endeq
where $\omega_1$ is the right half of $\omega$ and $c_{\omega_1}$ is its center.

Define
\begineq
\beta_l(x)=\sum_{\omega\in \mathcal{D}_l} \beta_{\omega}(x),
\endeq
note that 
\begineq
\beta_l(x+2^{-l})=\beta_l(x), \forall x\in [-2, 2-2^{-l}].
\endeq
Define 
\begineq
\gamma_l=\frac{1}{2}\int_{-1}^1 \beta_l(x+t)dt,
\endeq
because of the above periodicity, we know that $\gamma_l$ is constant for $x\in [-1,1]$, independent of $l$. Say $\gamma_l(x)=\delta>0$, hence 
\begineq
\frac{1}{\delta}\gamma_l(x) \mathds{1}_{[-1,1]}(x)=\mathds{1}_{[-1,1]}(x).
\endeq
Define another multiplier $\tilde{\beta}:\R\to \R$ with support in $[\frac{1}{2}, \frac{5}{2}]$ and $\tilde{\beta}(x)=1$ for $x\in [1,2]$. We define the corresponding multiplier on $\R^2$:
\begin{align*}
& \hat{m}_{k,\omega}(\xi_1, \xi_2)=\tilde{\beta}(2^{-k}\xi_2)\beta_{\omega}(\frac{\xi_1}{\xi_2})\\
& \hat{m}_{k,l,t}(\xi_1, \xi_2)=\tilde{\beta}(2^{-k}\xi_2)\beta_l(t+\frac{\xi_1}{\xi_2})\\
& \hat{m}_{k,l}(\xi_1, \xi_2)=\tilde{\beta}(2^{-k}\xi_2)\gamma_l(\frac{\xi_1}{\xi_2})
\end{align*}
Then what we need to bound can be written as 
\begin{align*}
\|\sum_{k\in \Z}\sum_{l\in \Z}H_l P_k(f)\|_p &=\|\int_{-1}^1 \sum_{k\in \Z}\sum_{l\ge 0} H_{k-l}(\frac{1}{\delta}m_{k,l}*f) dt\|_p\\
			&\le \int_{-1}^1 \|\sum_{k\in \Z} \sum_{l\ge 0} H_{k-l} (\frac{1}{\delta}m_{k,l,t}*f)\|_p dt,
\end{align*}
where the terms $H_l P_k$ for $l>k$ in the sum vanish as explained before.

So it suffices to prove a uniform bound on $t\in [-1, 1]$, w.l.o.g. we will just consider the case $t=0$, which is 
\begineq
\sum_{k\in \Z} \sum_{l\ge 0} H_{k-l}(m_{k,l,0}*f)=\sum_{k\in \Z} \sum_{l\ge 0} H_{k-l}([\tilde{\beta}(2^{-k}\xi_2)\beta_l(\frac{\xi_1}{\xi_2})]*f).
\endeq

{\bf Constructing the tiles:} For each $k\in \Z$ and $\omega\in \mathcal{D}_l$ with $l\ge 0$, let $\mathcal{U}_{k,\omega}$ be a partition of $\R^2$ by rectangles of width $2^{-k}$ and length $2^{-k+l}$, whose long side has slope $-c(\omega)$, where $c(\omega)$ is the center of the interval $\omega$. If $s\in \mathcal{U}_{k,\omega}$, we will write $\omega_s:=\omega$, and $\omega_{s,1}$ to be the right half of $\omega$, $\omega_{s, 2}$ the left half.

An element of $\mathcal{U}_{k, \omega}$ for some $\omega\in \mathcal{D}_l$ is called a ``tile''. Define $\varphi_{k,\omega}$ such that
\begineq
|\hat{\varphi}_{k,\omega}|^2=\hat{m}_{k,\omega},
\endeq
then $\varphi_{k,\omega}$ is smooth by our assumption on $\beta$ mentioned above.

For a tile $s\in \mathcal{U}_{k,\omega}$, define 
\begineq\label{wavelet}
\varphi_s(p):= \sqrt{|s|} \varphi_{k,\omega} (p-c(s)),
\endeq
where $c(s)$ is the center of $s$. Notice that 
\begineq
\|\varphi_s\|_2^2= \int_{\R^2} |s| \varphi_{k,\omega}^2=|s|\int_{\R^2} \hat{m}_{k,\omega}=1,
\endeq
i.e. $\varphi_s$ is $L^2$ normalized.

The construction of the tiles above by uncertainty principle is to localize the function further in space, for this purpose we need
\begin{lem}(\cite{Ba2})
\begineq
f*m_{k,\omega}(x)=\lim_{N\to \infty} \frac{1}{4 N^2} \int_{[-N,N]^2} \sum_{s\in \mathcal{U}_{k,\omega}} \langle f,\varphi_s(p+\cdot)\rangle \varphi_s(p+x)dp
\endeq
\end{lem}
\noindent The above lemma allows us to pass to the model sum
\begin{align*}
\sum_{k\in \Z} \sum_{l\ge 0} H_{k-l}(f* m_{k,l,0})= \sum_{k\in \Z}\sum_{l\ge 0}\sum_{\omega\in \mathcal{D}_l}\sum_{s\in \mathcal{U}_{k,\omega}}\langle f, \varphi_s\rangle H_{k-l}(\varphi_s),
\end{align*}
define 
\begineq
\psi_s=\psi_{-\log(length(s))},
\endeq
and 
\begineq\label{EE5.13}
\phi_s(x):= \int \check{\psi}_s(t)\varphi_s(x-tv(x))dt,
\endeq
then the model sum turns to 
\begineq\label{E5.14}
\sum_{k\in \Z}\sum_{l\ge 0}\sum_{\omega\in \mathcal{D}_l}\sum_{s\in \mathcal{U}_{k,\omega}}\langle f, \varphi_s\rangle \phi_s
\endeq

\begin{lem}\label{keylemma}
we have that $\phi_s(x)=0$ unless $-u(h(x))\in \omega_{s,2}$.
\end{lem}
The proof of the above lemma is by the Plancherel theorem, we just need to observe that the frequency support of $\psi_s$ and $\hat{\varphi}_s$ will be disjoint at the point $x$ unless $-u(h(x))\in \omega_{s,2}$.

\section{Boundedness of the commutator term: Proof of Proposition \ref{prop2.4}}\label{section6}

In this section we intend to prove that for any $p> 1$, there exists $\gamma_p>0$ such that
\begineq\label{EE5.1}
\left\|\sum_{k\in \Z} \left(H_{v, k-l} P_k(f)-\tilde{P}_k H_{v, k-l} P_k(f)\right)\right\|_p \lesim 2^{-\gamma_p l}\|f\|_p.
\endeq
If we expand the left hand side of the last expression to a model sum by the notations in Subsection \ref{subsection4.1}, \eqref{EE5.1} becomes
\begineq\label{EE5.2}
\left\|\sum_{k\in \Z}\sum_{\omega\in \mathcal{D}_l}\sum_{s\in \mathcal{U}_{k,\omega}}\langle f, \varphi_s\rangle (\phi_s-\tilde{P}_k \phi_s)\right\|_p \lesim 2^{-\gamma_p l} \|f\|_p.
\endeq
Observe that for a fixed point $x\in \R^2$, by Lemma \ref{keylemma}, the expression 
\begineq
\sum_{k\in \Z}\sum_{s\in \mathcal{U}_{k,\omega}}\langle f, \varphi_s\rangle (\phi_s-\tilde{P}_k \phi_s)(x)
\endeq 
can be non-zero for at most one $\omega\in \mathcal{D}_l$, which implies that
\begineq\label{FF5.4}
\begin{split}
& \left\|\sum_{k\in \Z}\sum_{\omega\in \mathcal{D}_l}\sum_{s\in \mathcal{U}_{k,\omega}}\langle f, \varphi_s\rangle (\phi_s-\tilde{P}_k \phi_s)\right\|_p\\
& \lesim \left(\sum_{\omega\in \mathcal{D}_l}\int_{\R^2} \left|\sum_{k\in \Z}\sum_{s\in \mathcal{U}_{k,\omega}}\langle f, \varphi_s\rangle (\phi_s-\tilde{P}_k \phi_s)\right|^p\right)^{1/p}
\end{split}
\endeq
From the right hand side of the above inequality, we see that \eqref{EE5.2} is reduced to separate $\omega\in \mathcal{D}_l$. Hence we just need to do the estimate for each $\omega$ separately. To be precise, we will prove 
\begin{lem}\label{lemma5.1}
Under the above notations, we have
\begineq\label{FF5.5}
\left\|\sum_{k\in \Z}\sum_{s\in \mathcal{U}_{k,\omega}}\langle f, \varphi_s\rangle (\phi_s-\tilde{P}_k \phi_s)\right\|_p \lesim 2^{- l}\|P_{\omega}f\|_p,
\endeq
where $P_{\omega}$ is the frequency projection operator given by
\begineq
\mathcal{F} P_{\omega}f (\xi_1, \xi_2)=\beta_{\omega}(\frac{\xi_1}{\xi_2}) \mathcal{F} f(\xi_1, \xi_2),
\endeq
and the constant in \eqref{FF5.5} is independent of $\omega\in \mathcal{D}_l$.
\end{lem}

\begin{lem}\label{lemma4.4}
We have the following bounds for the multiplier $\beta_{\omega}$:
\begineq\label{NN5.7}
\|P_{\omega}f\|_p \lesim \|f\|_p,
\endeq
for all $p\in (1, \infty)$, with the constant being independent of $\omega$.\\
\end{lem}

\noindent {\bf Finishing the proof of Proposition \ref{prop2.4}:} we substitute the estimates in Lemma \ref{lemma5.1} and Lemma \ref{lemma4.4} into \eqref{FF5.4} to obtain
\begineq\label{FF5.7}
\begin{split}
& \left\|\sum_{k\in \Z}\sum_{\omega\in \mathcal{D}_l}\sum_{s\in \mathcal{U}_{k,\omega}}\langle f, \varphi_s\rangle (\phi_s-\tilde{P}_k \phi_s)\right\|_p\\
& \lesim \left(\sum_{\omega\in \mathcal{D}_l} 2^{-pl}\|P_{\omega}f\|_p^p\right)^{1/p}\lesim 2^{-\frac{p-1}{p}\cdot l}\|f\|_p,
\end{split}
\endeq
which finishes the proof of Proposition \ref{prop2.4}.$\Box$

\begin{rem}
It has been proved by Demeter and Di Plinio in \cite{DD} that 
\begineq
\left(\sum_{\omega\in \mathcal{D}_l} \|P_{\omega}f\|_p^p\right)^{1/p}\lesim \|f\|_p,
\endeq
for $p\ge 2$, with the constant being independent of $l\in \N$. This will provide a better exponential decay in $l$ in the last inequality in \eqref{FF5.7}. However, here we do not need such orthogonality estimate but simply a triangle inequality.
\end{rem}

\subsection{Proof of Lemma \ref{lemma4.4}}

We first reduce the estimate to one single $\omega\in \mathcal{D}_l$ by applying the shearing transform. Suppose for the moment that we have proved \eqref{FF5.5} for $\omega=[0, 2^{-l}]$, by doing the following change of variables
\begineq\label{NN4.23}
x_1\to x_1, x_2\to x_2+\lambda x_1,
\endeq 
for the function $f$, the frequency variables are transformed into
\begineq
\xi_1\to \xi_1-\lambda \xi_2, \xi_2\to \xi_2.
\endeq 
This linear change of variables turns 
\begineq
P_{\omega'}f(\xi_1, \xi_2)=\mathcal{F}^{-1}\left( \beta_{\omega'}(\frac{\xi_1}{\xi_2}) \hat{f}(\xi_1, \xi_2) \right),
\endeq
which is the term on the left hand side of \eqref{NN5.7}, into
\begineq\label{NN4.26}
\mathcal{F}^{-1}\left( \beta_{\omega'}(\frac{\xi_1}{\xi_2}) \hat{f}(\xi_1-\lambda \xi_2, \xi_2) \right).
\endeq
If we denote
\begineq
\tilde{\xi}_1:=\xi_1-\lambda \xi_2, \tilde{\xi}_2:=\xi_2,
\endeq
the multiplier in \eqref{NN4.26} turns to
\begineq\label{NN4.28}
\beta_{\omega'}(\frac{\tilde{\xi}_1+\lambda\tilde{\xi}_2}{\tilde{\xi}_2})=\beta(2^{l+c}\frac{\tilde{\xi_1}}{\tilde{\xi}_2}+\lambda 2^{l+c}-2^{l+c}c_{\omega'_1}).
\endeq
So far it becomes clear that by taking $\lambda$ in \eqref{NN4.28} properly, we can apply the change of variables \eqref{NN4.23} to turn the projection operator $P_{\omega'}f$ for an arbitrary $\omega\in \mathcal{D}_l$ to $P_{\omega}f$, where $\omega=[0, 2^{-l}]$. \\

Next, we will reduce the estimate for all $l\in \N$ to the one simply for $l=0$. This can be done by applying the following unisotropic scaling symmetry:
\begineq
x_1\to \lambda x_1, x_2\to x_2,
\endeq
for the function $f$. Under the above change of variables, the Fourier transform of $f$ is transformed from $\hat{f}(\xi_1, \xi_2)$ to
\begineq
\frac{1}{\lambda}\hat{f}(\frac{\xi_1}{\lambda}, \xi_2).
\endeq
Correspondingly, the function $P_{\omega}f$ is changed to
\begineq
\begin{split}
& \int \beta_{\omega}(\frac{\xi_1}{\xi_2})\frac{1}{\lambda}\hat{f}(\frac{\xi_1}{\lambda}, \xi_2)e^{ix_1 \xi_1+ix_2 \xi_2}d\xi_1 d\xi_2\\
& =\int \beta_{\omega}(\frac{\lambda \xi_1}{\xi_2})\hat{f}(\xi_1, \xi_2)e^{i\lambda x_1 \xi_1+ix_2 \xi_2}d\xi_1 d\xi_2.
\end{split}
\endeq
Hence the multiplier $\beta_{\omega}(\xi_1/\xi_2)$ has the same $L^p$ norm with  $\beta_{\omega}(\lambda \xi_1/\xi_2)$. However, by the definition of $\beta_{\omega}$, we have
\begineq
\beta_{\omega}(\frac{\lambda \xi_1}{\xi_2})=\beta( \frac{2^{l+c}\lambda \xi_1}{\xi_2}-2^{l+c}c_{\omega_1}),
\endeq
which means that if we take $\lambda=2^{-l}$, the right hand side of the last expression becomes $\beta_{\omega_0}(\xi_1/\xi_2)$ where $\omega_0=[0, 1]$.\\

After the above inductions, we just need to prove \eqref{NN5.7} with $\omega_0=[0, 1]$. For $p=2$, the estimate is trivial due to Plancherel's theorem. For $p\neq 2$, if we denote by $P_k$ a Littlewood-Paley projection operator in the second variable, then by the Littlewood-Paley theory, we obtain
\begineq
\|P_{\omega_0}f\|_p \lesim \left\|\left(\sum_k \left|P_k P_{\omega_0}f\right|^2\right)^{1/2}\right\|_p.
\endeq
By the classical Calderon-Zygmund theory, it is not difficult to prove that 
\begineq
\left\|\left(\sum_k \left|P_k P_{\omega_0}f\right|^2\right)^{1/2}\right\|_{BMO}\lesim \|f\|_{\infty},
\endeq
and 
\begineq
\left\|\left(\sum_k \left|P_k P_{\omega_0}f\right|^2\right)^{1/2}\right\|_{1}\lesim \|f\|_{H^1}.
\endeq
Hence by interpolation, we obtain the desired estimate for all $p\in (1, \infty)$. So far we have finished the proof of Lemma \ref{lemma4.4}. $\Box$

\subsection{Proof of Lemma \ref{lemma5.1}}
By the same shearing transform as in \eqref{NN4.23}, we can reduce the estimate \eqref{FF5.5} for different $\omega$ to the one for a fixed $\omega$, say $\omega=[0, 2^{-l}]$. To prove \eqref{FF5.5}, by invoking duality, it is equivalent to prove that 
\begineq\label{EE5.3}
\int_{\R^2} \sum_{k\in \Z}\sum_{s\in \mathcal{U}_{k,\omega}}|\langle f, \varphi_s\rangle| \left|\left(\phi_s-\tilde{P}_k \phi_s\right) \cdot g\right| \lesim 2^{- l} \|f\|_p,
\endeq
where the function $g$ satisfies $\|g\|_{p'}\le 1$. By the Fubini theorem, the left hand side of \eqref{EE5.3} is equal to 
\begineq\label{EE5.10}
\begin{split}
& \int_{\R} \int_{\R^2}\sum_{k\in \Z}\sum_{s\in \mathcal{U}_{k,\omega}}|\langle f, \varphi_s\rangle| \left|\left(\phi_s(x)-\tilde{P}_k \phi_s(x)\right) \cdot g(x)\right| \delta(h(x)-t)dxdt\\
& = \int_{\R} \int_{\R^2}\sum_{k\in \Z}\sum_{s\in \mathcal{U}_{k,\omega}}|\langle f, \varphi_s\rangle| \left|\left(\phi_s(x)-\tilde{P}_k \phi_s(x)\right) \cdot g(x)\right| \delta(h_t(x)) \frac{|\nabla h_t(x)|}{|\nabla h(x)|}dxdt.
\end{split}
\endeq
By the bound on $\nabla h_t$ in \eqref{GG3.5} and our assumption on $\nabla h$ in the Main Theorem, the right hand side of \eqref{EE5.10} can be bounded by
\begineq\label{EE5.12}
\begin{split}
\int_{\R} \int_{\R^2}\sum_{k\in \Z}\sum_{s\in \mathcal{U}_{k,\omega}}|\langle f, \varphi_s\rangle| \left|\left(\phi_s(x)-\tilde{P}_k \phi_s(x)\right) \cdot g(x)\right| \delta(h_t(x)) dxdt.
\end{split}
\endeq
If we denote by $s_{m,n}$ the translation of the tile $s$ by $(m,n)$ units, which is
\begineq
s_{m,n}:=s-(m\cdot l_s, n\cdot w_s),
\endeq 
then the above \eqref{EE5.12} is equal to
\begineq\label{NN4.33}
\sum_{m,n}\int_{\R} \int_{\R^2}\sum_{k\in \Z}\sum_{s\in \mathcal{U}_{k,\omega}}|\langle f, \varphi_s\rangle| \mathbbm{1}_{s_{m,n}}(x)\left|\left(\phi_s(x)-\tilde{P}_k \phi_s(x)\right) \cdot g(x)\right| \delta(h_t(x)) dxdt.
\endeq 
By the notion of the adapted rectangles in Definition \ref{definition3.6}, we can replace $s_{m,n}$ by the slightly enlarged ``rectangle'' $\tilde{s}_{m,n}$ as from the definition it is clear that $\tilde{s}_{m,n}\supset s_{m,n}$. Moreover, in the following, we will only focus on the term $m=n=0$, as the other terms appear as the tail terms by the non-stationary phase method.\\

To proceed, we need the notion of Jones' beta numbers:
\begin{defi}
Fix a Lipschitz function $A:\R\to \R$. For each dyadic interval $I$, there exists a number $\alpha_I(A)\in \R$, such that if we denote
\begineq
\beta_{j_0}(I)=\sup_{x\in 3 j_0 I}\frac{|A(x)-A(c_I)-\alpha_I(A)(x-c_I)|}{|I|},
\endeq
where $j_0\in \N$ and $c_I$ denotes the center of $I$, then we will have the following Carleson type condition
\begineq
\sup_J \frac{1}{|J|}\sum_{I\subset J}\beta_{j_0}^2(I)|I|\lesim j_0^3 \|A\|_{Lip}^2.
\endeq
$\beta_{j_0}(I)$ and $\alpha_I(A)$ will be called the $j_0$-th beta number and the ``average slope'' for the Lipschitz function $A$ near the interval $I$ separately.
\end{defi}

The pointwise estimate in the following Lemma \ref{mainlemma2} will play a crucial role in the forthcoming calculation. To state this estimate, we need to make some preparations: for a fix $t\in \R$, we use the new coordinates system given by $(v_t, v_t^{\perp})$. For a tile $s$, we use $J(t, s)$ to denote the projection of $\Gamma_t \cap \tilde{s}$ on the new vertical axis $v^{\perp}_t$. Moreover for the interval $J(t, s)$, we let $J^D(t, s)$ denote one of the dyadic intervals (at most two) on the vertical axis such that 
\begineq
|J^D(t, s)|\in (8\cdot |J(t, s)|, 16\cdot |J(t, s)|]
\endeq
and
\begineq
|J^D(t, s)\cap J(t, s)|\ge |J(t, s)|/2.
\endeq
For the dyadic interval $J^D(t, s)$, we let $\Phi_{J^D(t, s)}$ denote the associated $L^2$ normalized Haar function.

\begin{lem}\label{mainlemma2}(\cite{Guo})
Fix $t\in \R$ and $s\in \mathcal{U}_{k,\omega}$ for some $\omega\in \mathcal{D}_l$, for $x\in \Gamma_t\cap \tilde{s}$, we have the pointwise estimate
\begineq\label{NN5.12}
|\phi_s(x)-\tilde{P}_k \phi_s(x)|\lesim \sum_{j_0\in \N}\frac{2^{-3l/2} 2^k \beta_{j_0}(J^D(t, s)) }{<j_0>^N},
\endeq
where $\beta_{j_0}(J^D(t, s))$ denotes the $j_0$-th beta number of $\Gamma_t$ near the dyadic interval $J^D(t, s)$.
\end{lem}

\begin{rem}
The proof of the above Lemma \ref{mainlemma2} in \cite{Guo} relies on those unnecessary parameters and auxiliary functions that we want to avoid by doing $\delta$-calculus. However, as we have promised in the introduction that we will carry out the whole argument in the language of $\delta$-calculus completely, we should also be able to prove Lemma \ref{mainlemma2} by doing so. This is postponed to the next subsection.
\end{rem}

Substitute the above estimate into the right hand side of \eqref{NN4.33} with $m=n=0$, we obtain
\begineq\label{EE5.15}
\sum_{j_0}\frac{2^{-l}}{<j_0>^N}\int_{\R} \int_{\R^2}\sum_{k\in \Z}\sum_{s\in \mathcal{U}_{k,\omega}}|\langle f, \varphi_s\rangle|2^k 2^{-l/2} \mathbbm{1}_{\tilde{s}}(x)\beta_{j_0}(J^D(t, s)) \cdot |g(x)| \delta(h_t(x)) dxdt.
\endeq
To proceed, we need the following
\begin{claim}\label{claim5.5}
Fix $t\in \R$, we have the following estimate
\begin{align*}
& \int_{\R^2} \left( \sum_{k\in \Z}\sum_{s\in \mathcal{U}_{k, \omega}}2^k 2^{-l/2}  |\langle f, \varphi_s\rangle| \mathbbm{1}_{\tilde{s}}(x)\beta_{j_0}(J^D(t, s)) \right) g(x)\delta(h_t(x))dx \\
& \lesim j_0^{3/2} \left(\int_{\R^2}\left(\sum_{k\in \Z}\sum_{s\in \mathcal{U}_{k, \omega}}|\langle f, \varphi_s\rangle|^2 \chi_s^2(x)\right)^{p/2}\delta(h_t(x))dx\right)^{1/p} \left(\int_{\R_2}|g(x)|^{p'}\delta(h_t(x))dx\right)^{1/p^{\prime}},
\end{align*}
where for $x=(x_1, x_2)$,
\begineq
\chi_s(x_1, x_2):=\frac{|s|^{-1/2}}{(1+(\frac{x_1-c_{s,1}}{l_s})^2+(\frac{x_2-c_{s,2}}{w_s})^2)^5},
\endeq
with $c_s=(c_{s, 1}, c_{s, 2})$ denoting the center of $s$, $l_s=2^{-k+l}$ the length and $w_s=2^{-k}$ the width. 
\end{claim}

We postpone the proof of the Claim \ref{claim5.5} till the end of this subsection and continue with the estimate of the term \eqref{EE5.15}. By Claim \eqref{claim5.5} and by applying H"older's inequality to $\int_{\R} dt$, the expression in \eqref{EE5.15} can be bounded by
\begineq\label{EE5.17}
\begin{split}
& \sum_{j_0}\frac{ j_0^{3/2}\cdot 2^{-l}}{<j_0>^N}\left(\int_{\R}\int_{\R^2}\left(\sum_{k\in \Z}\sum_{s\in \mathcal{U}_{k, \omega}}|\langle f, \varphi_s\rangle|^2 \chi_s^2(x)\right)^{p/2}\delta(h_t(x))dxdt\right)^{1/p}\\
& \lesim 2^{-l}\cdot \left\|\left(\sum_{k\in \Z}\sum_{s\in \mathcal{U}_{k, \omega}}|\langle f, \varphi_s\rangle|^2 \chi_s^2(x)\right)^{1/2}\right\|_p. 
\end{split}
\endeq
To bound the last expression, we need the following
\begin{lem}\label{lemma5.6}
We have the following variant of the square function estimate
\begineq\label{GG6.16}
\left\|\left(\sum_{k\in \Z}\sum_{s\in \mathcal{U}_{k, \omega}}|\langle f, \varphi_s\rangle|^2 \chi_s^2(x)\right)^{1/2}\right\|_p \lesim \|f\|_p.
\endeq
\end{lem}

\vspace{4mm}

\noindent {\bf Finishing the proof of Lemma \ref{lemma5.1}:} it is straightforward that, combined with \eqref{EE5.17}, Lemma \ref{lemma5.6} finishes the estimate of the expression \eqref{EE5.15}, thus the proof of Lemma \ref{lemma5.1}. $\Box$\\

\noindent {\bf Proof of Lemma \ref{lemma5.6}:} recall that in the estimate \eqref{GG6.16}, we have $\omega=[0, 2^{-l}]$. Now we want to reduce the estimate to the case $\omega_0=[0, 1]$ by applying the unisotropic scaling 
\begineq
x_1\to 2^l x_1, x_2\to x_2.
\endeq
Under the above change of variables, as has been explained in the proof of Lemma \ref{lemma4.4}, $\varphi_s$ for some $s\in \mathcal{U}_{k, \omega}$ is changed to $\varphi_{s'}$ for the corresponding $s'\in \mathcal{U}_{k, \omega_0}$ with $\omega_0=[0, 1]$. Moreover,  the function $\chi_s$ will also behave in the same way: 
\begineq\label{NN5.37}
\begin{split}
\chi_s(2^l x_1, x_2)& =\frac{|s|^{-1/2}}{(1+(\frac{2^l x_1-c_{s,1}}{l_s})^2+(\frac{x_2-c_{s,2}}{w_s})^2)^5}\\
& = \frac{|s|^{-1/2}}{(1+(\frac{x_1-2^{-l}c_{s,1}}{2^{-l}l_s})^2+(\frac{x_2-c_{s,2}}{w_s})^2)^5}.
\end{split}
\endeq
Recall that $l_s=2^l w_s$, hence the right hand side of \eqref{NN5.37} becomes a bump function with main support on a cube of side length $w_s$, which means that $\chi_s(2^l x_1, x_2)$ is equal to $\chi_{s'}$ for some $s'\in \mathcal{U}_{k, \omega_0}$  up to a normalization factor.\\

After the above reduction, we just need to prove \eqref{GG6.16} for $\omega=[0, 1]$. For the case $p=2$, by the orthogonality of the wavelet functions, we obtain
\begineq
\left(\int_{\R^2} \sum_{k\in \Z} \sum_{s\in \mathcal{U}_{k, \omega}} |\langle f, \varphi_s\rangle|^2 \chi_s^2\right)^{1/2} \lesim \|f\|_2.
\endeq
Moreover, by the classical Calderon-Zygmund theory, it is not difficult to prove the following endpoint estimates
\begineq
\left\|\left(\sum_{k\in \Z} \sum_{s\in \mathcal{U}_{k, \omega}} |\langle f, \varphi_s\rangle|^2 \chi_s^2\right)^{1/2}\right\|_{BMO} \lesim \|f\|_{\infty},
\endeq
and 
\begineq
\left\|\left(\sum_{k\in \Z} \sum_{s\in \mathcal{U}_{k, \omega}} |\langle f, \varphi_s\rangle|^2 \chi_s^2\right)^{1/2}\right\|_{1} \lesim \|f\|_{H^1}.
\endeq
Hence by interpolation, we can obtain all the expected $L^p$ estimate for \eqref{GG6.16} in the above Lemma \ref{lemma5.6}. $\Box$\\

\noindent {\bf Proof of Claim \ref{claim5.5}:} for a fixed $t\in \R$, for the summation on the left hand side of the estimate in Claim \ref{claim5.5}, we observe that
\begineq
\sum_{k\in \Z}\sum_{s\in \mathcal{U}_{k, \omega}}=\sum_{s:s\cap \Gamma_t\neq \emptyset},
\endeq 
as the term $\mathbbm{1}_{\tilde{s}}(x)$ will vanish  if $s\cap \Gamma_t=\emptyset$. We use the new coordinate system $(v_t, v^{\perp}_t)$, and write $x=x_1 v_t+x_2 v^{\perp}_t$, which will still be denoted as $x=(x_1, x_2)$ for the sake of simplicity. This turns the left hand side of the estimate in Claim \ref{claim5.5} into
\begineq\label{EE5.21}
\begin{split}
& \int_{\R}\int_{\R}\left( \sum_{s:s\cap \Gamma_t\neq \emptyset}|\langle f, \varphi_s\rangle| \mathbbm{1}_{\tilde{s}}(x_1, x_2)\beta_{j_0}(J^D(t, s)) 2^k 2^{-l/2} \right) g(x_1, x_2)\delta(h_t(x_1, x_2))dx_1dx_2\\
& = \sum_{s:s\cap \Gamma_t\neq \emptyset}2^k 2^{-l/2}  |\langle f, \varphi_s\rangle| \beta_{j_0}(J^D(t, s)) \int_{\R}\int_{\R}g(x_1, x_2)\mathbbm{1}_{\tilde{s}}(x_1, x_2)\delta(h_t(x_1, x_2))dx_1dx_2
\end{split}
\endeq

%

\noindent Notice that the integration on the right hand side of \eqref{EE5.21} can be estimated in the following way
\begin{align*}
& \left|\int_{\R}\int_{\R}g(x_1, x_2)\mathbbm{1}_{\tilde{s}}(x_1, x_2)\delta(h_t(x_1, x_2))dx_1dx_2\right| \\
& \lesim 2^{-k} \left[\int_{\R}g(x_1, \cdot)\delta(h_t(x_1, \cdot))dx_1\right]_{2J^D(t, s)},
\end{align*}
where for a function $G:\R\to \R$, $[G(\cdot)]_J$ denotes the average of the function $G$ on the interval $J\subset \R$.

Substitute the above bound into the right hand side of \eqref{EE5.21}, we obtain the following bound
\begin{align*}
\sum_{s:s\cap \Gamma_t\neq \emptyset}2^{-l/2} |\langle f, \varphi_s\rangle| \beta_{j_0}(J^D(t, s))\left[\int_{\R}g(x_1, \cdot)\delta(h_t(x_1, \cdot))dx_1\right]_{2J^D(t, s)}.
\end{align*}

To proceed, the idea is to view the above expression as a paraproduct. To do this, we need to find the right function such that it has the wavelet coefficient $2^{-l/2} |\langle f, \varphi_s\rangle| w_s^{-1/2}$, where $w_s=2^{-k}$ denotes the width of the tile $s$. This can be achieved by defining a function $F_t:\R\to \R$ such that
\begineq
F_t(x_2)=\sum_{s: s\cap \Gamma_t\neq \emptyset}2^{-l/2} w_s^{-1/2} \langle f, \varphi_s\rangle \Phi_{J^D(t, s)}(x_2),
\endeq 
where $\Phi_{J^D(t, s)}$ denotes the $L^2$ normalized Haar function associated to the dyadic interval $J^D(t, s)$.\\

By the $L^p$ boundedness of the paraproduct (see \cite{AHMTT} for example) and Jones' beta number condition that
\begineq
\sup_{s} \frac{1}{|J^D(t, s)|} \sum_{s': J^D(t, s')\subset J^D(t, s)}\beta_{j_0}^2(J^D(t, s')) w_s \lesim j_0^3,
\endeq
we obtain for any fixed $t\in \R$ that
\begineq
\begin{split}
& \sum_{s: s\cap \Gamma_t\neq \emptyset} 2^{-l/2} |\langle f, \varphi_s\rangle|  \beta_{j_0}(J^D(t, s))\left[\int_{\R}g(x_1, \cdot)\delta(h_t(x_1, \cdot))dx_1\right]_{2J^D(t, s)}\\
& =\sum_{s: s\cap \Gamma_t\neq \emptyset} 2^{-l/2} w_s^{-1/2} |\langle f, \varphi_s\rangle|  \beta_{j_0}(J^D(x, s)) w_s^{1/2} \left[\int_{\R}g(x_1, \cdot)\delta(h_t(x_1, \cdot))dx_1\right]_{2J^D(t, s)}\\
& \lesim j_0^{3/2} \|F_t(\cdot)\|_p \left\|\int_{\R}g(x_1, \cdot)\delta(h_t(x_1, \cdot))dx_1\right\|_{p'}\\
& \lesim j_0^{3/2} \|F_t(\cdot)\|_p \left(\int_{\R_2}|g(x)|^{p'}\delta(h_t(x))dx\right)^{1/p^{\prime}}\end{split}
\endeq

Hence what remains is to prove the following
\begin{claim}\label{claim5.7}
Under the above notations, we have
\begineq\label{EE5.28}
\|F_t(\cdot)\|_p \lesim \left(\int_{\R^2}\left(\sum_{k\in \Z}\sum_{s\in \mathcal{U}_{k, \omega}}|\langle f, \varphi_s\rangle|^2 \chi_s^2(x)\right)^{p/2}\delta(h_t(x))dx\right)^{1/p}.
\endeq
\end{claim}
\noindent {\bf Proof of Claim \ref{claim5.7}:} by the square function estimate, we obtain
\begineq\label{EE5.29}
\|F_t\|_p \lesim \left\|\left(\sum_{s:s\cap \Gamma_t\neq \emptyset}2^{-l}w_s^{-2} \langle f, \varphi_s\rangle^2 \mathbbm{1}_{J^D(t,s)}(\cdot)\right)^{1/2}\right\|_p.
\endeq
For the right hand side of \eqref{EE5.28}, again we use the new coordinate system $(v_t, v_t^{\perp})$ and denote $x=x_1 v_t + x_2 v_t^{\perp}$ as $x=(x_1, x_2)$ for the sake of simplicity.
Then the right hand side of \eqref{EE5.28} becomes
\begineq\label{EE5.30}
\begin{split}
& \left(\int_{\R^2}\left(\sum_{k\in \Z}\sum_{s\in \mathcal{U}_{k, \omega}}|\langle f, \varphi_s\rangle|^2 \chi_s^2(x_1, x_2)\right)^{p/2}\delta(h_t(x_1, x_2))dx_1dx_2\right)^{1/p}\\
& =\left(\int_{\R}\left(\int_{\R}\sum_{k\in \Z}\sum_{s\in \mathcal{U}_{k, \omega}}|\langle f, \varphi_s\rangle|^2 \chi_s^2(x_1, x_2)\delta(h_t(x_1, x_2))dx_1\right)^{p/2}dx_2\right)^{1/p}.
\end{split}
\endeq

If we compare the right hand side of \eqref{EE5.29} and \eqref{EE5.30}, we observe that the following pointwise estimate in $x_2$ will finish the proof of the claim: for any $x_2\in \R$ and any tile $s$ such that $s\cap \Gamma_t\neq \emptyset$, we have
\begineq
\begin{split}
 2^{-l}w_s^{-2} \langle f, \varphi_s\rangle^2 \mathbbm{1}_{J^D(t,s)}(x_2)\lesim \int_{\R}|\langle f, \varphi_s\rangle|^2 \chi_s^2(x_1, x_2)\delta(h_t(x_1, x_2))dx_1.
\end{split}
\endeq
But this follows easily from the definition of the function $\chi_s$.  Thus we have finished the proof of Claim \ref{claim5.7}. $\Box$

\subsection{Proof of Lemma \ref{mainlemma2}}


As we are fixing $t$ and trying to prove pointwise estimate for $x\in \Gamma_t$, we could always pretend that the vector field is constantly equal to $v_t$ on the whole plane. That is to say, if we define
\begineq\label{GG5.17}
\phi_s^t(x):=\int_{\R}\varphi_s(x-t v_t)\check{\psi}_{k-l}(t)dt, \forall x\in \R^2,
\endeq
we will have
\begineq
\phi_s^t(x)=\phi_s(x), \forall x\in \Gamma_t,
\endeq
and the advantage is that the vector field becomes the constant vector field $v_t$. In the following, we will stick to $\phi_s^t$ instead of $\phi_s$.\\

For a tile $s$ of dimension $w_s\times l_s$ with 
\begineq
l_s=2^l\cdot w_s,
\endeq 
for a point $x\in \Gamma_t\cap s$ with 
\begineq
v_t^{\perp}\in \omega_{s,2},
\endeq
we want to show that 
\begineq\label{BB6.55}
|\phi_s^t(x)-\tilde{P}_k \phi_s^t(x)|\lesim \sum_{j_0\in \N}\frac{2^{-3l/2}\cdot w_s^{-1} \beta_{j_0}(J^D(t, s)) }{<j_0>^N}.
\endeq
To proceed, we again turn to the new coordinate system $(v_t, v_t^{\perp})$, and write
\begineq
x\to x_1 v_t+x_2 v_t^{\perp}.
\endeq
By the definition of the operator $\tilde{P}_k$, the left hand side of \eqref{BB6.55} is equal to 
\begineq\label{GG5.22}
\phi_s^t(x_1, x_2)-\int_{\R}\left[\int_{\R}\phi_s^t(y_1, y_2)\delta(h_t(y_1, y_2))dy_1\right]\psi_k(x_2-y_2)dy_2.
\endeq
We approximate $\Gamma_t\cap s$ by the line of the ``average slope'' in the definition of Jones' $\beta$-number, and call it $l_{s,t}$. Moreover, we define another auxiliary function $L_{s,t}$ associated to the line $l_{s,t}$ in a similar way to $h_t$:
\begineq
\text{If for some $y\in \Gamma_t$ we have $z-y=d\cdot v_t$, then we set $L_{s,t}(z)=d$.}
\endeq
The crucial observation is that 
\begineq
\begin{split}
& \int_{\R}\left[\int_{\R}\phi_s^t(y_1, y_2)\delta(L_{s,t}(y_1, y_2))dy_1\right]\psi_k(x_2-y_2)dy_2\\
& =\int_{\R}\phi_s^t(y_1, x_2)\delta(L_{s,t}(y_1, x_2))dy_1.
\end{split}
\endeq
Substitute the above identity into \eqref{GG5.22} to obtain
\begineq\label{GG5.25}
\begin{split}
& \phi_s^t(x_1, x_2)-\int_{\R}\phi_s^t(y_1, x_2)\delta(L_{s,t}(y_1, x_2))dy_1...\\
& ... - \int_{\R}\left[\int_{\R}\phi_s^t(y_1, y_2)\left(\delta(h_t(y_1, y_2))-\delta(L_{s,t}(y_1, y_2))\right)dy_1\right]\psi_k(x_2-y_2)dy_2. 
\end{split}
\endeq
Notice that for $x=(x_1, x_2)\in \Gamma_t$, we have 
\begineq
\phi_s^t(x_1, x_2)=\int_{\R}\phi_s^t(y_1, x_2)\delta(h_t(y_1, x_2))dy_1,
\endeq
by substituting which into \eqref{GG5.25} we obtain
\begineq\label{GG5.27}
\begin{split}
& \int_{\R}\phi_s^t(y_1, x_2)\left[\delta(h_t(y_1, x_2))-\delta(L_{s,t}(y_1, x_2))\right]dy_1...\\
& ...  - \int_{\R}\left[\int_{\R}\phi_s^t(y_1, y_2)\left(\delta(h_t(y_1, y_2))-\delta(L_{s,t}(y_1, y_2))\right)dy_1\right]\psi_k(x_2-y_2)dy_2. 
\end{split}
\endeq
Observe that the latter term in the above expression is just a Littlewood-Paley projection of the former term, hence it should be expected that these two terms can be handled in a similar way. In \cite{Guo} it is indeed shown to be this case, hence in the following we will focus on the former term of \eqref{GG5.27}.

By the definition of $\phi_s^t$ in \eqref{GG5.17}, we obtain
\begineq\label{GG5.28}
\begin{split}
& \int_{\R}\phi_s^t(y_1, x_2)\left[\delta(h_t(y_1, x_2))-\delta(L_{s,t}(y_1, x_2))\right]dy_1\\
&=\int_{\R}\int_{\R}\varphi_s(y_1-t, x_2)\check{\psi}_{k-l}(t)dt\left[\delta(h_t(y_1, x_2))-\delta(L_{s,t}(y_1, x_2))\right]dy_1.
\end{split}
\endeq
If we denote
\begineq
d:=h_t(y_1, x_2)-L_{s,t}(y_1, x_2),
\endeq
then the right hand side of \eqref{GG5.28} turns to 
\begineq
\begin{split}
& \int_{\R}\int_{\R}\left(\varphi_s(y_1-t, x_2)-\varphi_s(y_1+d-t, x_2)\right)\check{\psi}_{k-l}(t)dt \delta(h_t(y_1, x_2))dy_1\\
& =\int_{\R}\int_{\R}\varphi_s(y_1-t, x_2)\left(\check{\psi}_{k-l}(t)-\check{\psi}_{k-l}(t+d)\right)dt \delta(h_t(y_1, x_2))dy_1.
\end{split}
\endeq
Hence by the definition of Jones'  beta numbers that 
\begineq
|d|\lesim w_s\cdot \beta_0(J^D(t, s)),
\endeq
and by applying the fundamental theorem to $\check{\psi}_{k-l}$, we conclude the desired estimate in Lemma \eqref{mainlemma2}. $\Box$

Shaoming Guo, Institute of Mathematics, University of Bonn\\
\indent Address: Endenicher Allee 60, 53115, Bonn\\
\indent Email: shaoming@math.uni-bonn.de

\end{document}